\documentstyle[12pt]{article}

\renewcommand{\theequation}{\arabic{section}.\arabic{equation}}

\def\be{\begin{equation}}
\def\ee{\end{equation}}
\def\bea{\begin{eqnarray}}
\def\eea{\end{eqnarray}}

\def\1{\'{\i}}                           
\def\R{{\rm I\kern-.2em R}}

\begin{document}

\thispagestyle{empty}

\ 
\vspace{3cm}

\begin{center} {\LARGE{\bf{Quantum Heisenberg--Weyl Algebras}}} 
\end{center}

\bigskip\bigskip

\begin{center} Angel Ballesteros$^\dagger$, Francisco J.
Herranz$^\dagger$ and Preeti Parashar$^\ddagger$
\end{center}

\begin{center} {\it {  $^\dagger$ Departamento de F\1sica, Universidad
de Burgos} \\   Pza. Misael Ba\~nuelos, 
E-09001-Burgos, Spain}
\end{center}

\begin{center} {\it {  $^\ddagger$ SISSA, Via Beirut 2-4, 
34014 Trieste, Italy}}
\end{center}

\bigskip\bigskip\bigskip

\begin{abstract} 
All Lie bialgebra structures on the Heisenberg--Weyl
algebra $[A_+,A_-]=M$ are classified and explicitly quantized. The
complete list of quantum Heisenberg--Weyl algebras so obtained includes
new multiparameter deformations, most
of them being of the non-coboundary type. 
\end{abstract} 

\newpage

\setcounter{equation}{0}

\renewcommand{\theequation}{\arabic{equation}}


 A Hopf algebra deformation of a universal enveloping algebra $U g$
defines in a unique way a Lie bialgebra structure $(g,\delta)$ on $g$
 \cite{PC}. The cocommutator $\delta$ provides the first order terms
in the deformation of the coproduct, and can be seen as the natural
tool to classify quantum algebras.  Moreover, this well known statement
suggests the relevance of the inverse  problem, i.e., to find a method
to construct, given an arbitrary Lie bialgebra, a Hopf algebra
quantization of it. 

This question has been addressed recently in
\cite{osc}, where a very general  construction of a deformed
coassociative coproduct linked to a given Lie bialgebra $(g,\delta)$
has been presented. Such Lie bialgebra quantization formalism,
inspired by the paper
\cite{LM} (see also \cite{Lyak}),  has been shown to be universal for
the oscillator algebra:  multiparametric coproducts corresponding to
all coboundary oscillator Lie bialgebra structures
can be obtained in that way (for  the oscillator algebra non-coboundary
structures do not exist
\cite{oscGos}). To complete  the structure of quantum algebras,
deformed commutation rules can be found by imposing the homomorphism
condition for the coproduct (counit and antipode can be also easily
derived).

In this letter we show that all  Heisenberg--Weyl Lie bialgebras can
be completely quantized by making use of this formalism. This result
enhances the advantages of such an approach in order to obtain a full
chart of Hopf algebra deformations of physically relevant algebras. 

Firstly, we find the most general form  of all families
of Heisenberg--Weyl Lie bialgebras. It is remarkable that, in contrast to
the oscillator case, now there exists only one coboundary bialgebra
among them. Afterwards, it is shown  how all these  Lie bialgebras can
be classified and ``exponentiated" to get the quantum coproducts by
means of the formalism introduced in
\cite{osc}. We also find all deformed  commutation rules, thus
obtaining a complete list of quantum deformations of this algebra,
whose properties are briefly commented. This exhaustive description is
fully complementary with respect to the quantum group results already
known either from a Poisson--Lie construction
\cite{Kuper} or from an
$R$-matrix approach \cite{HLR}. 

Let us fix the notation. The Heisenberg--Weyl Lie algebra $h_3$ is
generated by $A_+$, $A_-$ and $M$ with Lie brackets
\be
 [A_-,A_+]=M,\qquad 
[M,\cdot\,]=0 .
\label{aa}
\ee
A $3\times 3$ real matrix representation $D$ of (\ref{aa}) is given by:
\be 
 D(A_+)=\left(\begin{array}{ccc}
 0 &0 & 0 \\ 0 & 0 & 1 \\ 0 & 0 & 0 
\end{array}\right),\quad
 D(A_-)=\left(\begin{array}{ccc}
 0 &1 & 0 \\ 0 & 0 & 0 \\ 0 & 0 & 0 
\end{array}\right),\quad 
D(M)=\left(\begin{array}{ccc}
 0 &0 & 1 \\ 0 & 0 & 0 \\ 0 & 0 & 0 
\end{array}\right).
\label{ab}
\ee 
The expression for a generic element of  the Heisenberg--Weyl group
$H_3$ coming from this representation is: 
\be 
 D(T)= \exp\{m D(M)\}\exp\{a_- D(A_-)\}
\exp\{a_+ D(A_+)\}  =\left(\begin{array}{ccc}
 1 &a_-   & m  + a_- a_+ \\ 0 & 1 & a_+ \\ 0 & 0 & 1 
\end{array}\right) ,
\label{ac}
\ee 
and the group law for the coordinates $m$, $a_-$ and $a_+$  is
obtained by means of matrix multiplication $D({T}'')=D({T}')\cdot
D({T})$:
\be 
 m''=m+m' 
- a_- a'_+  ,\quad 
 a''_+=a'_+ + a_+  ,\quad 
 a''_-=a'_- + a_- .
\label{ad}
\ee


Heisenberg--Weyl Lie bialgebras $(h_3,\delta)$ will be defined by the
cocommutator
$\delta:h_3\to h_3\otimes h_3$ such that

\noindent 
i) $\delta$ is a 1--cocycle, i.e.,
\be
\delta([X,Y])=[\delta(X),\, 1\otimes Y+ Y\otimes 1] + 
[1\otimes X+ X\otimes 1,\, \delta(Y)], \quad \forall \,X,Y\in
h_3. \label{ba}
\ee

\noindent 
ii) The dual map $\delta^\ast:h_3^\ast\otimes h_3^\ast \to
h_3^\ast$ is a Lie bracket on $h_3^\ast$.

From ii), we consider an arbitrary skewsymmetric cocommutator:
\bea
&&\delta(A_-)=a_1\, A_-\wedge A_+ + 
a_2\, A_-\wedge M + a_3\, A_+ \wedge
M,\cr 
&&\delta(A_+)=b_1\, A_-\wedge A_+ + b_2\, A_-\wedge M + b_3\, A_+
\wedge M,\cr 
&&\delta(M)=c_1\, A_-\wedge A_+ + c_2\, A_-\wedge M + c_3\, A_+
\wedge M,
\label{bc}
\eea
where $a_i,b_i,c_i$ ($i=1,2,3$) are real parameters.
If we impose on (\ref{bc}) the cocycle condition (\ref{ba}) we obtain:
\be
c_1=0,\qquad c_2=b_1,\qquad c_3=-a_1 .
\label{bd}
\ee
Since the dual $h^\ast_3$ with generators  $\{m ,  a_-,a_+\}$  must be
a  Lie algebra, the Jacobi identity on the bracket $\delta^\ast$ gives
rise to two additional conditions:
\be
a_1(b_3 - a_2) - 2 b_1 a_3 = 0,\qquad 
b_1(a_2 - b_3) - 2 a_1 b_2 = 0.
\label{bf}
\ee
Hence, the most general Heisenberg--Weyl bialgebra has commutation
relations (\ref{aa}) and cocommutators
\bea
&&\delta(A_-)=a_1\, A_-\wedge A_+ 
+ a_2\, A_-\wedge M + a_3\, A_+ \wedge
M,\cr 
&&\delta(A_+)=b_1\, A_-\wedge A_+ + b_2\, A_-\wedge M + b_3\, A_+
\wedge M,\cr 
&&\delta(M)=  b_1\, A_-\wedge M - a_1\, A_+
\wedge M,
\label{bg}
\eea
with the six parameters $a_i$, $b_i$ verifying (\ref{bf}).

It is also known that the dual Lie bracket
$\delta^\ast$ gives the linear part of the (unique) Poisson--Lie
structure on the group linked to
$\delta$ \cite{Dr}. Therefore, starting from the
classification of Poisson--Lie Heisenberg groups given in \cite{Kuper}
and taking into account the change of local coordinates on the
Heisenberg group
\be
x_1=a_-,\qquad x_2=a_+,\qquad x_3=m + a_-\,a_+,
\label{bhb}
\ee
it is straightforward to prove that the full Poisson--Lie
bracket associated to $\delta$ reads
\bea
&&\{a_-,a_+\}=a_1\, a_- + b_1\, a_+,\cr
&&\{a_-,m\}=
a_2\, a_- + b_2\, a_+ + b_1\, m - \frac {a_1}{2}\, a_-^2,\cr
&&\{a_+,m\}=a_3\, a_- + b_3 \,a_+ - a_1\, m + \frac {b_1}{2}\, a_+^2 .
\label{bh}
\eea
In other words, if (\ref{ad}) is read as a coproduct on Fun($H_3$), it
is easy to check that the group law turns out to be a Poisson algebra
homomorphism with respect to (\ref{bh}).

Finally, let us find out for which
values of the parameters we have coboundary Lie bialgebras. So, we
investigate the most general  skewsymmetric element $r$ of $h_3\otimes
h_3$ such that 
\be
\delta(X):=[1\otimes X + X \otimes 1,\,  r],\quad 
X\in h_3, \label{bi}
\ee
defines a Lie
bialgebra. This is equivalent  to imposing the Schouten bracket
$[[r,r]]$ to be a solution of the modified classical Yang--Baxter
equation (YBE)
\be
[X\otimes 1\otimes 1 + 1\otimes X\otimes 1 +
1\otimes 1\otimes X,[[r,r]]\, ]=0, \quad X\in h_3.
\label{bj}
\ee
Explicitly, we consider three real-valued coefficients
$\xi,\beta_+,\beta_-$ and write:
\be 
 r= \xi\, A_+\wedge A_- + \beta_+ \, A_+\wedge M
 + \beta_-  A_- \wedge M.
\label{bl}
\ee
The Schouten bracket of this element is given by
 \be 
  [[r,r]]= - \xi^2 \,M\wedge A_+\wedge A_-.
\label{bm}
\ee 
This bracket is found to fulfill automatically the modified classical
YBE (\ref{bj}). Therefore, (\ref{bl}) is always a classical $r$-matrix.
The cocommutator (\ref{bi}) derived from it reads
\be
\delta(A_+)= - \xi\, A_+\wedge M,\qquad 
\delta(A_-)= - \xi\, A_-\wedge M,\qquad 
\delta(M)= 0 .
\label{bn}
\ee
Thus, we conclude that there exists only one non-trivial coboundary
Heisenberg--Weyl Lie bialgebra which is characterized by
\be
a_1=a_3=b_1=b_2=0,\qquad a_2=b_3=- \xi.
\label{bo}
\ee
The case $\xi=0$
gives rise to a solution of the classical YBE, but now the 
cocommutator vanishes.


Let us go back to the four-parameter family of bialgebras given by
(\ref{bf}) and (\ref{bg}). It is easy to check that equations (\ref{bf})
have three disjoint types of solutions:

\noindent  Type I$_+$: $a_1\neq 0$, $b_2=-\,a_3\,b_1^2/a_1^2$,
$b_3=a_2+2\,b_1\,a_3/a_1$ and
$a_2,a_3,b_1$ arbitrary. 

\noindent  Type I$_-$: $a_1=0$,  $b_1\neq 0$, $a_3=0$, $a_2=b_3$ and
$b_2,b_3$ arbitrary. 

\noindent  Type II: $a_1=0$, $b_1=0$ and $a_2,a_3,b_2,b_3$ arbitrary.

So, we have three (multiparametric) families of Lie bialgebras. To
quantize them, we have to check that, within each family \cite{osc}:

\noindent a)  There exists some set $\{H_i\}$ of commuting
generators of $g$ such that $\delta(H_i)=0$  (these will be the
primitive generators after quantization).

\noindent b) For the remaining generators $X_j$, their
cocommutator $\delta(X_j)$ must only contain  terms of the form
$X\wedge H$ (neither $X_l\wedge X_m$ nor
$H_n\wedge H_p$ contributions are allowed).

Finally, we have to take into account that two Lie bialgebra
structures of a Lie algebra $g$ are equivalent if there exists an
automorphism of $g$ that transforms one into the other. As we shall
see, some automorphisms of the Heisenberg algebra will help us to get
bialgebras fulfilling conditions a) and b).


\noindent $\bullet$ {\bf Type I$_+$:} This is a
family of Lie bialgebras which has, for general  values of the
parameters, no primitive generator
$\delta(H)=0$. However, if we define
\be
A_+':=A_+ - \frac{b_1}{a_1}\,A_- + \left(\frac{b_1\,a_3}{a_1^2} 
+ \frac{a_2}{a_1} \right)\,M,\qquad a_1\neq 0,
\label{qa}
\ee
it is immediate to check that, in this new basis, the Type I$_+$
bialgebras have the following cocommutator:
\bea
&& \delta(A_-)=- a_1\, A_+'\wedge A_-  + a_3\, A_+' \wedge
M,\cr
&& \delta(A_+')=0,\cr
&& \delta(M)=   - a_1\, A_+' \wedge M .
\label{cb}
\eea
The automorphism (\ref{qa}) has shown the parameters $b_1$ 
and $a_2$ to be superfluous. 

The coproduct that quantizes the resultant biparametric family
(\ref{cb}) can be now obtained: firstly, we see that this family of
bialgebras verify conditions a) and b) with
$A_+'$ being the primitive generator (from now on, we shall write
$A_+$ instead of $A_+'$). Following \cite{osc} we write the
non-vanishing cocommutators in (\ref{cb}) in the matrix form:
\be
 \delta\left(\begin{array}{c}
A_- \\ M  
\end{array}\right)=
\left(\begin{array}{cc}
-a_1 A_+ & a_3 A_+   \\ 0 & -a_1 A_+
\end{array}\right)\dot\wedge \left(\begin{array}{c}
A_- \\ M   
\end{array}\right). 
\label{cc}
\ee
In this way, the coproduct for non-primitive  generators will be
formally given by:
\be  
 \Delta\left(\begin{array}{c}
A_- \\ M  
\end{array}\right)=
\left(\begin{array}{cc}
1 &0   \\ 0 & 1 
\end{array}\right)
\dot\otimes \left(\begin{array}{c}
A_- \\ M  
\end{array}\right)   +
\sigma\left( \exp\left\{\left(\begin{array}{cc}
 a_1 A_+ & -a_3 A_+   \\ 0 &  a_1 A_+
\end{array}\right)\right\}\dot\otimes \left(\begin{array}{c}
A_- \\ M    
\end{array}\right)\right) .
\label{cd}
\ee
By computing explicitly the exponential, we find that
\bea
&&\Delta(A_+)=1\otimes A_+ + A_+ \otimes 1,\qquad 
 \Delta(M)=1\otimes M +M \otimes e^{a_1A_+},\cr
&&\Delta(A_-)=1\otimes A_- + A_- \otimes e^{a_1A_+} -
a_3 M \otimes A_+\, e^{a_1A_+}. 
\label{ed}
\eea
The next step is the search for deformed commutation rules compatible
with (\ref{ed}). They turn out to be
\be 
 [A_-,A_+]=M,\qquad 
[A_-,M]=\frac {a_1}2 M^2 ,\qquad 
[A_+,M]=0 .
\label{eg}
\ee 
Finally, counit and antipode are deduced
\be
\epsilon(X)=0,\qquad X\in\{A_-,A_+,M\},
\label{ee}
\ee
\bea
&&\gamma(A_+)=-A_+,\qquad \gamma(M)=-M \,e^{-a_1A_+},\cr 
&&\gamma(A_-)=-A_- \, e^{-a_1A_+} - a_3 M\, A_+ \, e^{-a_1A_+} , 
\label{ef}
\eea
and the Hopf algebra $U_{a_1,a_3}(h_3)$ that quantizes
the family of (non-coboundary) Heisenberg--Weyl  bialgebras (\ref{cb})
is obtained.

It is remarkable that in this quantum deformation the parameter $a_3$
is not involved in the deformed  commutation rules. Recall that
(\ref{eg}) was firstly obtained in \cite{HLR} starting from a quantum
Heisenberg group and by
applying a duality method (coproduct (\ref{ed}) could not be
found). 

On the other hand, a
physically suggestive observation  comes from the fact that the
generator
$M$ is neither central nor primitive (recall the role that the
non-primitive mass generator of  quantum extended Galilei algebra
plays in one-dimensional magnon systems \cite{magn}). The central
element
$\cal C$ is now
\be
{\cal C}=M\,e^{-a_1\,A_+/2}.
\label{eeg}
\ee
This element labels the following differential realization of
$U_{a_1,a_3}(h_3)$:
\be
A_+=x,\qquad A_-=\lambda\,e^{a_1\,x/2}\,\partial_x,\qquad
M=\lambda\,e^{a_1\,x/2},
\label{eeh}
\ee
where $\lambda$ is the eigenvalue of $\cal C$. Note also that, by
introducing $\cal C$ as a new generator instead of $M$, relations
(\ref{eg}) turn into 
\be 
 [A_-,A_+]={\cal C}\,e^{a_1\,A_+/2},\qquad 
[A_-,{\cal C}]=0 ,\qquad 
[A_+,{\cal C}]=0 .
\label{eei}
\ee 

Finally note that, if ${\tilde A}$,  ${\tilde A}_+$ and ${\tilde A}_-$
are the generators of the non-standard quantum deformation $U_z
sl(2,\R)$
\cite{sl}, the quantum Heisenberg algebra $U_{a_1,0}(h_3)$ can be
obtained as the contraction $\varepsilon \to 0$ defined by
\be
M=-\varepsilon\,{\tilde A},\qquad A_+={\tilde A}_+, 
\qquad A_-=\varepsilon\,{\tilde A}_-,\qquad a_1=2\,z.
\label{eej}
\ee
As it could be expected from the non-coboundary character of
$U_{a_1,0}(h_3)$, the universal $R$-matrix of $U_z sl(2,\R)$ diverges
under (\ref{eej}).
 

\noindent $\bullet$ {\bf Type I$_-$:} After  specializing the
corresponding parameters we find a three-parameter cocommutator also
with no primitive generators. However, the definition of $A_-'$ by
means of the automorphism
\be
A_-':=A_- - \frac{b_3}{b_1}\,M ,\qquad b_1\neq 0,
\label{qb}
\ee
implies that this family of Lie bialgebras is given by
\bea
&&\delta(A_-')= 0 ,\cr 
&&\delta(A_+)=b_1\, A_-'\wedge A_+ + b_2\, A_-'\wedge M ,\cr 
&&\delta(M)=  b_1\, A_-'\wedge M .
\label{bgh}
\eea
In particular, the parameter $b_3$ has been  reabsorbed, and
(\ref{bgh}) can be quantized. Moreover, this Type I$_-$ structures are
essentially the same as the Type I$_+$ (\ref{cb}), but reversing the
role of $A_-$ and
$A_+$. Once again, another Heisenberg algebra
automorphism given by
\be
A_+\to A_-,\qquad A_-\to A_+,\qquad M\to -M,
\label{ca}
\ee
would make both types of bialgebras  explicitly equivalent. Therefore,
we omit the explicit quantization leading to the algebra
$U_{b_1,b_2}(h_3)$.


\noindent $\bullet$ {\bf Type II:} If $a_1$ and $b_1$ vanish,
the cocommutator (\ref{bg}) reads:
\bea
&& \delta(A_-)=  a_2\, A_-\wedge M + a_3\, A_+ \wedge M,\cr
&& \delta(A_+)=  b_2\, A_-\wedge M + b_3\, A_+ \wedge M,\cr 
&& \delta(M)=0 .
\label{eh}
\eea
In this case, $M$ is the primitive generator  and no extra
manipulation is needed in order to quantize this family of bialgebras.
We write (\ref{eh}) in matrix form:
\be
 \delta\left(\begin{array}{c}
A_- \\ A_+  
\end{array}\right)=
\left(\begin{array}{cc}
-a_2 M & -a_3 M   \\ -b_2 M & -b_3 M
\end{array}\right)\dot\wedge \left(\begin{array}{c}
A_- \\ A_+ 
\end{array}\right), 
\label{ei}
\ee
hence, the corresponding coproduct is given by:
\be  
 \Delta\left(\begin{array}{c}
A_- \\ A_+
\end{array}\right)=
\left(\begin{array}{cc}
1 &0   \\ 0 & 1 
\end{array}\right)
\dot\otimes \left(\begin{array}{c}
A_- \\ A_+ 
\end{array}\right)   +
\sigma\left( \exp\left\{\left(\begin{array}{cc}
a_2 M & a_3 M   \\ b_2 M & b_3 M
\end{array}\right)\right\}\dot\otimes \left(\begin{array}{c}
A_- \\ A_+  
\end{array}\right)\right) .
\label{ej}
\ee

Although the four parameters describing this quantum algebra are
arbitrary, in order to derive the commutation rules compatible with
(\ref{ej}) it will suffice to write
\be
E:=\exp\left\{\left(\begin{array}{cc}
a_2 M & a_3 M   \\ b_2 M & b_3 M
\end{array}\right)\right\}=\left(\begin{array}{cc}
E_{11}(M) & E_{12}(M)   \\ E_{21}(M) & E_{22}(M)
\end{array}\right).
\label{ejb}
\ee
In this way, the explicit quantum coproduct will be
\bea 
&&  \Delta(M)=1\otimes M +M \otimes 1,\nonumber\\
&&  \Delta(A_-)=1\otimes A_- + A_- \otimes E_{11}(M) + A_+ \otimes
E_{12}(M),\nonumber\\
&&  \Delta(A_+)=1\otimes A_+ + A_+ \otimes E_{22}(M) + A_- \otimes
E_{21}(M).
\label{ejc}
\eea
Now, by taking into account the following property
\be
E_{11}(M)\,E_{22}(M)-E_{12}(M)\,E_{21}(M)=e^{(a_2+b_3)\,M}
\ee
it is straightforward to prove that the four-parameter coproduct
(\ref{ejc}) is an algebra homomorphism with respect to the deformed
commutation rules
\be 
 [A_-,A_+]=\frac { e^{(a_2+b_3)\,M}-1}{a_2+b_3} ,\qquad 
[A_-,M]=0 ,\qquad 
[A_+,M]=0 .
\label{en}
\ee
Due to the preservation of $M$ as central element, counit and antipode
are easily deduced. These operations complete the obtention of the
multiparametric quantum algebra
$U_{a_2,a_3,b_2,b_3}(h_3)$. These Type II  quantizations were studied
in
\cite{Lyak} with no reference to Lie bialgebra structures.

The well-known coboundary quantization is a particular subcase
with $a_3=b_2=0$ and $a_2=b_3=-\xi$. A  universal $R$-matrix (which is
not a solution of the quantum YBE) for it was obtained
in \cite{BCH}.

\bigskip

\noindent {\bf Acknowledgements}

\medskip

A.B. and F.J.H. are
partially supported by DGICYT (Project PB94-1115) from the
Ministerio de Educaci\'on y Ciencia de Espa\~na.



\begin{thebibliography}{40}

\bibitem{PC} Pressley A and Chari V 1990 {\it Nucl. Phys. B (Proc.
Suppl.)} {\bf 18A}  207

\bibitem{osc}
    Ballesteros A and  Herranz F J  1996
{\it J. Phys. A: Math. Gen.} {\bf 29}  4307

\bibitem{LM}
     Lyakhovsky V and    Mudrov A 1992
{\it J. Phys. A: Math. Gen.} {\bf 25}  L1139

\bibitem{Lyak}
     Lyakhovsky V 1994
{\it Zapiski Nauchn. Semin.} POMI V.209; 1995 
{\it Group-like structures in quantum-Lie algebras  and the procedure
of quantization}, in ``Quantum Groups, Formalism and Applications", J.
Lukierski {\it et al}, Eds., p 93, Polish Scientific Publishers.

\bibitem{oscGos}
    Ballesteros A and  Herranz F J  1996
{\it Quantum Oscillator Hopf Algebras}, communication  to the XXI
ICGMTP,  Goslar.

\bibitem{Kuper} Kuperschmidt B A 1993 
{\it J. Phys. A: Math. Gen.} {\bf 26} L929

\bibitem{HLR}
     Hussin V, Lauzon A and   Rideau G 1994
{\it Lett. Math. Phys.} {\bf 31}  159

\bibitem{Dr} 
Drinfel'd V G 1983
{\it Sov. Math. Dokl.} {\bf 27}  68


\bibitem{magn}
Bonechi F, Celeghini E,   Giachetti R,   Sorace E and   
Tarlini M  1992 {\it Phys. Rev B} 
{\bf 46} 5727 

\bibitem{sl} 
       Ballesteros A and Herranz F J 1996
{\it J. Phys. A: Math. Gen.} {\bf 29}  L311

\bibitem{BCH} 
       Ballesteros A, Celeghini E, Herranz F J,   
del Olmo M A and   Santander M 1994
{\it J. Phys. A: Math. Gen.} {\bf 27}  L369

\end{thebibliography}
\end{document}